**Using the Cover Guessing Game Strategy to Predict the Results of Stochastic Processes**

By James D. Stein (California State University, Long Beach)

**Abstract** – Cover originally proposed a game in which Bob wrote two numbers on pieces of paper and allowed Alice to pick one and look at it. Cover then showed – as Blackwell is thought to have shown previously – that Alice has a strategy involving the selection of a random number which enables her to guess correctly with a probability greater than ½ which of the two numbers was larger. Central to Cover's reasoning was that Alice could have picked either number to look at with equal probability.

In this paper, we investigate the situation in which Bob obtains the two numbers as the result of an experiment based on a probability distribution, such as the spin of a roulette wheel. We then modify this to envision Alice and Bob conducting separate experiments based on probability distributions. Alice conducts one trial of her experiment and uses Cover's strategy to guess whether the result from her experiment will be larger or smaller than the result from Bob's experiment – whenever he chooses to conduct a single trial of it. The results of this paper include

(1) If the two numbers are observations of the same experiment, Alice merely has to observe one trial and apply the Cover strategy to have a probability greater than ½ of guessing the larger value – even if the second trial has not yet been conducted

(2 There are many non-trivial examples where Alice and Bob conduct different experiments, and Alice can observe one trial of her experiment and apply the Cover strategy to have a probability greater than ½ of guessing the larger value – even if many of the specifics of the other experiment have not yet been determined and a trial of it has not yet been conducted.

**Introduction**

There is some controversy over the origin of the Cover Guessing Game ([2]), Portnoy ([4]) refers to two papers of David Blackwell published in the 1950s as the origin of the basic idea of using a random variable to improve one's guessing probability in certain situations. Cover stated the problem in its present form, and Wapner ([5]) gave it the name Blackwell's Bet, by which it is more widely known.

However, no one seems to have asked, "Where did Bob get these two numbers?" We shall show that if Bob's numbers are the result of two consecutive trials of a probabilistic experiment, such as the spin of a roulette wheel or the determination of the half-life of a radioactive atom in a jar of radioactive atoms of different species, Alice need only see the result of the first of the two trials to guess the larger with probability greater than ½.

Cover's original guessing game and the relevant strategy are simple. Two envelopes each contain a different number. Pick one envelope at random (it is assumed that each envelope has an equal probability of being picked), observe the number inside it, and pick a random number from any distribution non-zero on all open intervals. If the random number is smaller than the observed number, guess that the observed number is the larger of the two, and conversely.

We call the random number the pointer, because if one draws a vector from the observed number to the pointer, the vector points in the direction of what the guessing strategy posits as the unobserved number. If the pointer is less than the observed number, the vector points to what the guessing strategy thinks is the smaller number, and if the pointer is larger than the observed number, the vector points to what the guessing strategy thinks is the larger number.

The proof of the viability of this strategy is well-known but is included for completeness. Let S be the smaller of the two numbers in the envelopes, and L the larger. Let p be the probability that the random number is less than S, and q the probability that the random number is larger than L. With probability ½, the observed number will be S, and the probability that we will make the correct guess is 1-p. With probability ½, the observed number will be L, and the probability that we will make the correct guess is 1-q. The combined probability of making the correct guess is ½ (1-p) + ½ (1-q) = ½ + ½ (1 – (p+q)). This is ½ + ½ the probability that the pointer falls in the gap between S and L, and is thus larger than ½ as long as the pointer has a non-zero probability of lying in any open interval.

The difference between this situation, which is studied by Gnedin ([4]) using the tools of game theory, and the one we examine in this paper is that one of the two numbers could be handed to Alice without her having any chance to have picked the other number – including the possibility that the other number has not yet been determined.

. Alice's strategy is specified in advance. One can imagine that Alice is performing Experiment A here and now, and observing the value the experiment delivers, whereas Bob is on the other side of the Universe and won't conduct his experiment for billions of years. Nonetheless, there are many situations – as indicated in the Abstract – in which Alice can conduct her experiment and use the pointer to guess with probability greater than ½ whose experimental value is larger.

Throughout this paper, the phrase 'Experiment A is larger' is to be interpreted as the result of the single trial of Experiment A is larger than the single trial of Experiment B.

It is also somewhat surprising that the experiments need have nothing in common. Alice could be examining the spin of a roulette wheel, whereas Bob could be determining the half-life of a radioactive atom. As Pythagoras said, all is number.

**Section I – Identical Experiments**

It does not surprise us that if Alice and Bob were both given jars of 100 balls numbered 1 through 100, Alice could correctly guess with probability greater than ½ whose ball has the larger number. Alice simply looks at the number on the ball she has picked, and compares it with 50.5, the median of the distribution. If it is less than 50.5, she guesses her ball is the smaller, and conversely.

What may be surprising is that the pointer method works no matter what the distribution of balls in the jars are – even if nothing whatever is known about the distribution. This can actually be looked at as an extension of Cover's original guessing game. In that game, there were two trials of the identical experiment – that's how Bob got the numbers he wrote on separate pieces of paper. Alice was allowed to look at the result obtained from either one of the trials. In this case, Alice is only allowed to look at the result of the first trial, but is nonetheless able to correctly guess which is larger with a probability greater than ½.

**Theorem 1 –** Suppose Alice and Bob are performing identical experiments with outcomes $v_1 < \ldots < v_n$, where the outcome $v_i$ has probability $q_i$. Alice performs her experiment first, obtains a result, and then uses a pointer from a distribution that is non-zero on open sets. Alice guesses that her result is greater than Bob's if the pointer is less than her result and that Bob's result is greater if the pointer is greater than her result. We assume that if Alice and Bob have the same result, the probability that Alice guesses successfully is ½. Then P(SG), the probability that Alice will guess successfully, is greater than 1/2.

**Proof** - Let p be the pointer, and suppose $p < v_1$. Then Alice always guesses the result of Experiment A is larger. When she and Bob obtain the same result, which happens with probability $\sum_{i=1}^{n} {q_i}^2$ , Alice guesses correctly with probability ½. If Alice obtains result $v_i$, she guesses correctly if Bob obtains result $v_j$ where $j < i$, this happens with probability $q_i(q_1 + \ldots + q_{i-1})$. So Alice's probability of guessing successfully is

$$\frac{1}{2}\sum_{i=1}^{n} {q_i}^2 + q_2q_1 + q_3(q_1 + q_2) + \cdots + q_n(q_1 + \cdots + q_{n-1})$$

$$= \frac{1}{2}(q_1 + \cdots + q_n)^2 = \frac{1}{2}$$

Now suppose $p > v_n$. Then Alice always guesses the result of Experiment B is larger. When she and Bob obtain the same result, which happens with probability $\sum_{i=1}^{n} {q_i}^2$ , Alice guesses correctly with probability ½. If Alice obtains result $v_i$, she guesses correctly if Bob obtains result $v_j$ where $j > i$, this happens with probability $q_i(q_{i+1} + \ldots + q_n)$. So Alice's probability of guessing successfully is

$$\frac{1}{2}\sum_{i=1}^{n} {q_i}^2 + q_{n-1}q_n + q_{n-2}(q_{n-1} + q_n) + \cdots + q_1(q_2 + \cdots + q_n)$$

$$= \frac{1}{2}(q_1 + \cdots + q_n)^2 = \frac{1}{2}$$

If $v_i < p < v_{i+1}$, assume that the two outcomes of the two trials of the experiment are $v_j$ and $v_k$, with $j < k$. It is equally likely that the result of Experiment A is $v_j$ and the result of Experiment B is $v_k$ as it is that the result of Experiment A is $v_k$ and the result of Experiment B is $v_j$.

If $k \le i$

| Experiment A result | Experiment B result | Alice's guess | Guess result |
|---|---|---|---|
| $v_j$ | $v_k$ | B larger | right |
| $v_k$ | $v_j$ | B larger | wrong |

The guess is correct with probability ½.

If $j > i$

| Experiment A result | Experiment B result | Alice's guess | Guess result |
|---|---|---|---|
| $v_j$ | $v_k$ | A larger | right |
| $v_k$ | $v_j$ | A larger | wrong |

The guess is correct with probability ½.

If $j \le i$, $k > i$

| Experiment A result | Experiment B result | Alice's guess | Guess result |
|---|---|---|---|
| $v_j$ | $v_k$ | B larger | right |
| $v_k$ | $v_j$ | A larger | right |

The guess is always correct. Since this case happens with non-zero probability, the theorem is proved. █

We would obtain the same result, that the probability of a successful guess is greater than ½, if we simply ignored ties rather than assigning them a value of ½. The computation is essentially the same.

The technique in Theorem 1 also suffices to prove the result in the case where each experiment has a countable number of possible results, regardless of whether the results are bounded, bounded above or below, or unbounded.

**Ex. 1 (Coins and Balls) –** Suppose Alice and Bob each have two jars with identical contents. The first jar contains an unknown number of coins with unknown heads probabilities, the second jar contains an unknown number of numbered balls with unknown integers upon them. Alice picks one coin from the first jar and one ball from the second jar. She then flips the coin she picked the number of times indicated on the ball she picked, and records the number of heads. She then uses a pointer with a non-zero probability on each open interval to make her guess as to which number is larger – her number of heads, or Bob's, when he later performs the same experiment. By Theorem 1, her probability of guessing correctly is greater than ½.

**Ex. 2 (Half-Lives)** – Alice has a jar filled with radioactive atoms with assorted half-lives of different species. Assume there are a large number of atoms of each species. Alice picks an atom and determines its half-life. She can then use Cover's strategy to guess with probability greater than ½ whether the second atom she picks has a longer or shorter half-life.

The conclusion of Theorem 1 also holds for bounded continuous distributions.

**Theorem 2** – Assume the results of Experiments A and B (which are the same experiment) are selected from a distribution on [c,d] with density function $\rho(x) > 0$ on [c,d]. The probability that a result lies between x and x+dx is $\rho(x)dx$. Let p be a pointer taken from a distribution that is non-zero on open sets. Let x be the result of Experiment A, y the result of Experiment B. Alice is allowed to look at x but not y, and compare x to p. She guesses that $x < y$ if $p > x$, and guesses that $x > y$ if $p < x$.

(1) If $p < c$ or $p > d$, the probability of a successful guess is 0.5.

(2) If $c < p < d$, the probability of a successful guess is $> \frac{1}{2}$.

**Proof –** If $p < c$, x is always greater than p, and so Alice always guesses that $x > y$, which will be correct with probability 0.5. A similar argument holds if $p > d$.

Now assume $c < p < d$, and consider Fig. 1 as being in the xy-plane.

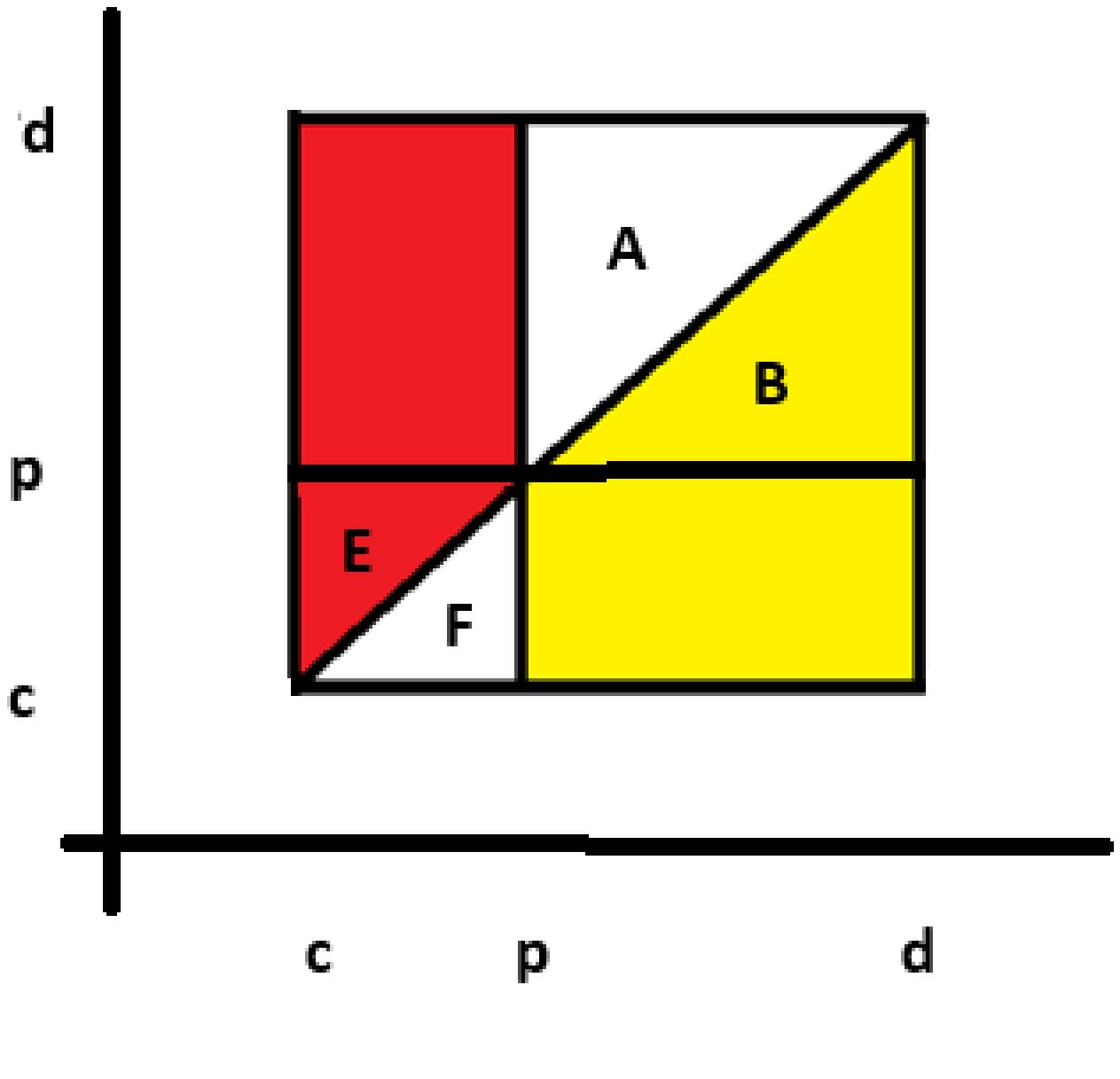

**Fig. 1**

If $x < p$, then guess that y is the larger of x and y. The probability of a successful guess is

$$\int_c^p \int_x^d \rho(y)\rho(x)dydx$$

If $x > p$, then guess that x is the larger of x and y. The probability of a successful guess is

$$\int_p^d \int_c^x \rho(y)\rho(x)dydx$$

The probability of a successful guess is therefore

$$\int_c^p \int_x^d \rho(y)\rho(x)dydx + \int_p^d \int_c^x \rho(y)\rho(x)dydx$$

The first integral is over the region colored red, the second integral is over the region colored yellow, so the integrand is taken over the entire colored area. The symmetry of the integrand ρ(x)ρ(y) ensures that the integral over triangle E equals the integral over triangle F, and the integral over triangle A equals the integral over triangle B. Since the integrand is non-negative and the colored area is greater than the white area, the integral over the colored area is greater than the integral over the colored area. Since the sum of the integral over the colored area and the integral over the white area is 1, the result follows. █

It seems likely that the boundedness condition can be dropped.

**Pointer Optimization**

Obviously, Alice is interested in improving her probability of a successful guess. As this section shows, if Alice can duplicate Bob's experiments she can ensure that P(SG) > ½ even if she knows nothing about the experiment other than that she has the same one that Bob does. But if she knows the full description of the experiment – both probabilities and outcomes – she can often maximize P(SG) through judicious selection of the pointer.

**Ex. 2**

| A values | 2 | 4 | 6 | 14 | 18 | | | | | |
|---|---|---|---|---|---|---|---|---|---|---|
| A probs | 0.4 | 0.3 | 0.1 | 0.1 | 0.1 | | | | | |
| | | | | | | | | | | |
| B values | 2 | 4 | 6 | 14 | 18 | | | | | |
| B probs | 0.4 | 0.3 | 0.1 | 0.1 | 0.1 | | | | | |
| | | | | | | | | | | |
| Pointers | 1 | 3 | 5 | 7 | 9 | 11 | 13 | 15 | 17 | 19 |
| P(SG) | 0.5000 | 0.8333 | 0.7917 | 0.7222 | 0.7222 | 0.7222 | 0.7222 | 0.6250 | 0.6250 | 0.5000 |
| | | | | | | | | | | |
| Uniform | Pointer | Mean | | 0.6604 | | | | | | |

The top row consists of the values of outcomes of Experiment A; directly below it are the probabilities of these outcomes. The next two rows consist of similar data for Experiment B (in this case, the two experiments are identical).

The possible pointers that Alice could decide to use are the odd numbers from 1 through 19; below each pointer is the probability of a successful guess from using that pointer. The uniform pointer mean is simply the average value of the probabilities of a successful guess for all 10 pointers.

Given total knowledge of the experiment, Alice would obviously choose to stick with 3 as a pointer, ensuring the correct prediction of the larger outcome 5 times of every 6 that the result does not end in a tie. If Alice is only aware that the outcomes lie between 1 and 19, the uniform random pointer selected as described above still selects the larger number approximately 2 out of every 3 times that a tie does not result.

## Mean vs Median

Although the mean attracts most of the attention in probability and statistics, the median may well play the greater role in pointer optimization, as the following example demonstrates.

**Ex. 3A**

| | | | | | | | |
|---|---|---|---|---|---|---|---|
| A values | 1 | 2 | 3 | 18 | | Mean | 6 |
| A probs | 0.25 | 0.25 | 0.25 | 0.25 | | | |
| | | | | | | | |
| B values | 1 | 2 | 3 | 18 | | Mean | 6 |
| B probs | 0.25 | 0.25 | 0.25 | 0.25 | | | |
| | | | | | | | |
| Pointers | 0.5 | 1.5 | 2.5 | 3.5 | 4.5 | 17.5 | 18.5 |
| P(SG) | 0.5000 | 0.7500 | 0.8333 | 0.7500 | 0.7500 | 0.7500 | 0.5000 |

The mean of the experiment is 6, but the median is 2.5, and for this example the optimal pointer is located at the median rather than the mean.

There is a natural tendency to focus on the median for pointer selection  The reason for this is straightforward – at least for uniform distributions such as in Ex. 3A..  When the pointer is located between the values obtained in the trials of Experiments A and B, Alice will always make the correct guess.  When the pointer is located at the median, Alice will make the correct guess whenever the two experimental results fall on different sides of the median – as well as in approximately half the cases where the two values are on the same side of the median. This also shows up in the Cover Paradox, as was noted in the Introduction.

However, it is possible to skew the distribution to favor the mean.  The following slight modification of Ex. 3A shows the crossover point for this example/

**Ex. 3B**

| A values | 1 | 2 | 3 | 18 | | Mean | 8.4 |
|---|---|---|---|---|---|---|---|
| A probs | 0.2 | 0.2 | 0.2 | 0.4 | | | |
| | | | | | | | |
| B values | 1 | 2 | 3 | 18 | | Mean | 8.4 |
| B probs | 0.2 | 0.2 | 0.2 | 0.4 | | | |
| | | | | | | | |
| Pointers | 0.5 | 1.5 | 2.5 | 3.5 | 4.5 | 17.5 | 18.5 |
| P(SG) | 0.5 | 0.722222 | 0.833333 | 0.833333 | 0.833333 | 0.833333 | 0.5 |

**When There Isn't Full Knowledge of the Experiment**

By "full knowledge", we mean knowledge of what the outcomes may be, and what their probabilities are. It is of interest (at least, to the author) to determine how to optimize the pointer in the absence of full knowledge.

Possibly the least one might know about the experiment is the range of values of the outcome. Actually, even if one does not know the exact range, one could assign lower and upper bounds beyond the range, and try a pointer which samples this range in some fashion – such as a uniform or normal distribution. We include calculations later for some situations in which we can compute the exact results.

**Section II – Experiments with Interwoven Outcomes**

In the previous section, we assumed that Alice and Bob were conducting the same experiment. In this section, we assume that the values of the two experiments are interwoven; some of the outcomes of Experiment A are between outcomes of Experiment B, and vice-versa.

**When Alice Is "Drawing Dead"**

There are situations in which no matter what pointer Alice tries, it is impossible for her to obtain P(SG) > ½, as the following example shows.

**Ex. 4**

| | | | | | | | | | | |
|---|---|---|---|---|---|---|---|---|---|---|
| **A values** | **2** | **6** | **10** | **11** | **18** | | | | | |
| **A probs** | **0** | **0.4** | **0** | **0.6** | **0** | | | | | |
| | | | | | | | | | | |
| **B values** | **3** | **5** | **8** | **12** | **16** | | | | | |
| **B probs** | **0.3** | **0.2** | **0** | **0.1** | **0.4** | | | | | |
| | | | | | | | | | | |
| **Pointers** | 1 | 3 | 5 | 7 | 9 | 11 | 13 | 15 | 17 | 19 |
| **P(SG)** | **0.5000** | **0.5000** | **0.5000** | **0.5000** | **0.5000** | **0.5000** | **0.5000** | **0.5000** | **0.5000** | **0.5000** |
| | | | | | | | | | | |
| **Uniform** | **Pointer** | **Mean** | | **0.5000** | | | | | | |

The probabilities associated with the values of Experiment B are symmetric with respect to each of the values of Experiment A. For each value of Experiment A, the probability is 0.5 that the value of Experiment B will be smaller, and the probability is 0.5 that the value of Experiment B will be larger.

The following two examples not only demonstrate this, but bring up an interesting point.

**Ex. 5**

| A values | 1 | 5 | 9 | 13 | 17 | | | | | | |
|---|---|---|---|---|---|---|---|---|---|---|---|
| A probs | 0 | 0.2 | 0.3 | 0.3 | 0.2 | | | | | | |
| | | | | | | | | | | | |
| B values | 3 | 7 | 11 | 15 | 19 | | | | | | |
| B probs | 0.1 | 0.3 | 0.2 | 0.3 | 0.1 | | | | | | |
| | | | | | | | | | | | |
| Pointers | 0 | 2 | 4 | 6 | 8 | 10 | 12 | 14 | 16 | 18 | 20 |
| P(SG) | 0.5000 | 0.5000 | 0.5000 | 0.6600 | 0.6600 | 0.7200 | 0.7200 | 0.6600 | 0.6600 | 0.5000 | 0.5000 |
| | | | | | | | | | | | |
| Uniform | Pointer | Mean | | 0.5982 | | | | | | | |

Contrast this with the following example.

**Ex. 6**

| A values | 1 | 6 | 9 | 14 | 17 | | | | | | |
|---|---|---|---|---|---|---|---|---|---|---|---|
| A probs | 0 | 0.3 | 0.24 | 0.26 | 0.2 | | | | | | |
| | | | | | | | | | | | |
| B values | 3 | 7 | 11 | 15 | 19 | | | | | | |
| B probs | 0.1 | 0.3 | 0.2 | 0.3 | 0.1 | | | | | | |
| | | | | | | | | | | | |
| Pointers | 0 | 2 | 4 | 6 | 8 | 10 | 12 | 14 | 16 | 18 | 20 |
| P(SG) | 0.4620 | 0.4620 | 0.4620 | 0.7020 | 0.7020 | 0.7500 | 0.7500 | 0.6980 | 0.6980 | 0.5380 | 0.5380 |
| | | | | | | | | | | | |
| Uniform | Pointer | Mean | | 0.6147 | | | | | | | |

Assume that Alice is aware of all the information in the above two examples, but is only allowed to use a uniformly selected random pointer from among the pointer choices. She places an even-money bet on the outcome of her guess.

Should Alice prefer the situation in Ex. 5 or Ex. 6? This is the classic 'minimax vs expectation' dilemma. Ex. 4 has the larger expectation – but Alice could have been unlucky enough to select either 0, 2 or 4 as a pointer and incur a long-term average loss. Ex. 5 has a smaller expectation, but the worst that would happen to Alice long-term is that she breaks even.

Notice that in both examples, the sum of P(SG) with a pointer of 0 and P(SG) with a pointer of 20 is 1. There is a simple reason for this; with a pointer of 0, Alice will always guess that the outcome of Experiment A is the larger, and with a pointer of 20 she will always guess that the outcome of Experiment B is the larger. One of these guesses must be correct every time the game is played, as the two experiments have no outcomes in common.

Notice also that P(SG) is the same for any pointer between the two adjacent values of Experiment A, as any two pointers satisfying this restriction will result in the same collection of values of Experiment A above or below the pointer – and it is those collections that determine P(SG).

We need to decide how to adjudicate Alice's guess when the outcome of Experiments A and B are the same. In some instance, it will be useful to adopt the Las Vegas approach that ties are considered as "no bet"; in some instances it will be useful to affix a win probability of ½ to Alice's guess in this case (as was done in the proof of Theorem 1). Whichever of these we adopt will not affect whether P(SG) > ½, although it will affect the actual value of P(SG). We are generally interested in whether P(SG) > ½, although there will be instances when we will compute the actual value. In those cases, we will specify which alternative we have adopted.

**When Any Choice of Pointer Yields P(SG) $\geq$ ½**

**Def.** We say that Experiment A strictly predicts Experiment B if every choice of pointer results in P(SG) $\geq$ ½ and for some choice of pointer, P(SG) > ½.

**Theorem 3** – (a) For every Experiment A with a finite number of outcomes, there is an Experiment B with outcomes having non-zero probabilities between every pair of outcomes from Experiment A that is strictly predicted by Experiment A.

(b) For every Experiment B with a finite number of outcomes, there is an Experiment A with outcomes having non-zero probabilities between every pair of outcomes from Experiment B there is an Experiment A with interwoven outcomes that strictly predicts Experiment B.

**Proof** – (a) Assume $A = \{v_k: k = 1,2,\ldots,N\}$, $v_1 < v_2 < \ldots < v_N$ and $P(v_k) = a_k$. Choose $w_1 < v_1$, choose $w_k$ such that $v_{k-1} < w_k < v_k$ for $k=1,2,\ldots,N$, and choose $w_{N+1} > v_N$. We first show we can define $b_k = P(w_k)$ such that for any pointer $< v_1$, $P(SG) = ½$.

Observe that the guessing strategy mandates that we always guess the value delivered by Experiment A is the larger. Given that Experiment A delivers the value $v_k$, the strategy guesses correctly if Experiment B delivers any of the values $w_1, \ldots, w_k$. So

$$P(SG) = a_1b_1 + a_2(b_1 + b_2) + \ldots + a_N(b_1 + \ldots + b_N)$$

$$= (a_1 + a_2 + \ldots + a_N)b_1 + (a_2 + a_3 + \ldots + a_N)b_2 + \ldots + a_Nb_N = ½ \qquad (*)$$

We need to ensure that there is a solution with $b_1 + \ldots + b_N < 1$. This can be written as

$b_1 + (1-a_1)b_2 + \ldots + (1 – (a_1+a_1+\ldots+a_{N-1}))b_N = ½$. For $2 \le k \le N$, choose $b_k = 1/(2N)$

This ensures that $b_2 + \ldots + b_N = (N-1)/(2N) < ½$, so $(1-a_1)b_2 + \ldots + (1 – (a_1+a_1+\ldots+a_{N-1}))b_N < ½$ as well. We can now choose $b_1 = ½ - (1-a_1)b_2 + \ldots + (1 – (a_1+a_1+\ldots+a_{N-1}))b_N$. So $0 < b_1 < ½$, and consequently $b_1 + (b_2 \ldots + b_N) < ½ + ½ = 1$.Now choose $b_{N+1} = 1 – (b_1 + \ldots + b_N)$ to complete the proof. Notice that, by the remark following Ex. 2, this also ensures that $P(SG) = ½$ for any pointer $> v_N$.

Finally, assume that the pointer p satisfies $v_{k-1} < p < v_k$. If the outcome of Experiment A is $v_j$ with $j < k$, Alice guesses that Experiment B has the larger value, otherwise she guesses that Experiment A has the larger value. So

$$P(SG) = a_1(b_2 + \ldots + b_{N+1}) + \ldots + a_{k-1}(b_k + \ldots + b_{N+1}) + a_k(b_1 + \ldots + b_k) + \ldots + a_N(b_2 + \ldots + b_N)$$

If p increases so that $v_k < p < v_{k+1}$, the net change to P(SG) is that if the outcome of Experiment A is $v_k$, Alice now guesses that Experiment A has the larger outcome. The net change to P(SG) is $a_k(b_{k+1} + \ldots + b_N) - a_k(b_1 + \ldots + b_k) = a_k(b_{k+1} + \ldots + b_N - (b_1 + \ldots + b_k))$. As k increases, the factor of $a_k$ is initially positive and results in increasing P(SG), but as fewer terms are added and more subtracted. It eventually becomes negative and P(SG) begins to decline.

(b) We employ essentially the same strategy as in the proof of (a). Define $\{v_k: k = 1,2,\ldots,N+1\}$ such that $v_1 < w_1$, $v_{N+1} > w_N$, and $w_k < v_{k+1} < w_{k+1}$ for $k = 1,2,\ldots,N-1$. Let p be a pointer with $p < v_1$. We want to define $P(v_k) = a_k$ such that P(SG) = ½. Notice that $p < v_k$ for all k, so we always guess that Experiment A has the larger value. Consequently, $P(SG) = a_2b_1 + a_3(b_1 + b_2) + \ldots + a_{N+1}(b_1 + b_2 + \ldots b_N) = a_2b_1 + a_3(b_1 + b_2) + \ldots + a_{N+1}$. We need to choose $a_2, \ldots ,a_{N+1}$ to ensure that this sum = ½, and that $a_1 + \ldots + a_{N+1} = 1$.

Choose $a_k = 1/(2N)$ for $k=2,3,\ldots,N$. Then $a_2 + \ldots + a_N = (N-1)/(2N) < ½$, so $a_2b_1 + a_3(b_1 + b_2) + \ldots + a_N(b_1 + b_2 + \ldots + b_{N-1}) < ½$. Let $a_{N+1} = ½ - (a_2b_1 + a_3(b_1 + b_2) + \ldots +$

$a_N(b_1 + b_2 + \ldots + b_{N-1}))$. Then P(SG) = ½ and$( a_2 + \ldots + a_N) + a_{N+1} < (N-1)/(2N) + ½ < 1$. So choose $a_1 = 1 – (a_2 + \ldots + a_{N+1})$ to conclude that P(SG) = ½ in this case. As in the proof of (a), P(SG) = ½ for any pointer $> v_{N+1}$, and the proof that P(SG) increases and then decreases as the pointer increases is similar to that of (a).. █

The next Theorem shows that experiments in which the probabilities display a weakened form of mirror symmetry can strictly predict experiments if the outcomes of Experiments A and B are interwoven.

**Theorem 4** – Suppose that $w_1 < v_1 < \ldots < w_n < v_n < v_{n+1} < w_{n+1} < \ldots < v_{2n} < w_{2n}$, where the $v_k$ represent the outcomes of Experiment A and the $w_k$ represent the outcomes of Experiment B. Suppose further that $P(v_k) = a_k = a_{2n+1-k}$ and $P(w_k) = b_k = b_{2n+1-k}$ (the probabilities are a weakened form of mirror-image symmetric). Then Experiment A strictly predicts Experiment B.

**Proof –** The arguments given in Theorem 3 show that it is only necessary to show that **i**f the pointer p is less than $v_1$, then $P(SG) = ½$. The decision principle guesses that the outcome from Experiment A is always the larger of the two outcomes, as the pointer is always less than the outcome of Experiment A. The probability P(SG) is given by

$$P(SG) = a_1b_1 + a_2(b_1 + b_2) + \ldots + a_n(b_1 + \ldots + b_n) + a_n(b_1 + \ldots + b_n)$$

$$+ a_{n-1}(b_1 + \quad + b_n + b_n) + \ldots + a_1(b_1 + \quad + b_n + b_n + \ldots + b_2)$$

Regroup to add the first term to the last, the second term to the next-to-last, etc. This gives

$$P(SG) = a_1(b_1 + \quad + b_n + b_n + \ldots + b_1) + \ldots + a_n(b_1 + \quad + b_n + b_n + \ldots + b_1) = a_1 + \ldots + a_n = ½.$$

If the pointer is greater than $v_{2n}$, then the decision principle guesses that the outcome from Experiment B is always the larger of the two outcomes, as the pointer is always greater than the outcome of Experiment A. If the above computation is thought of as going from left to right along the real line for k increasing from 1 to 2n, it can also be thought of as going from right to left along the real line for k decreasing from 2n to 1. So $P(SG) = ½$. The proof of Theorem 3 shows that for some pointer between $v_1$ and $v_{2n}$, we will have $P(SG) > ½$. █

**Section III – Uniform Distributions**

In this section, we will assume that all distributions – the two Experiments and the pointer – are uniformly distributed on intervals. This enables a relatively straightforward calculation of P(SG).

**Lemma 1** – Assume that the outcomes from Experiments A and B are uniform distributions on [0,1]. Let the pointer $p \in R$. Let x and y be the outcomes of Experiments A and B respectively..

(1) If $p < 0$ or $p > 1$, the probability of a successful guess is 0.5.

(2) If $0 < p < 1$, the probability of a successful guess is $½ + (p - p^2)$.

**Proof –** (1) If $p < 0$, x is always greater than p, and so we will always guess that $x > y$, which will be correct with probability 0.5. A similar argument holds if $p > 1$.

(2) Consider the following diagram.

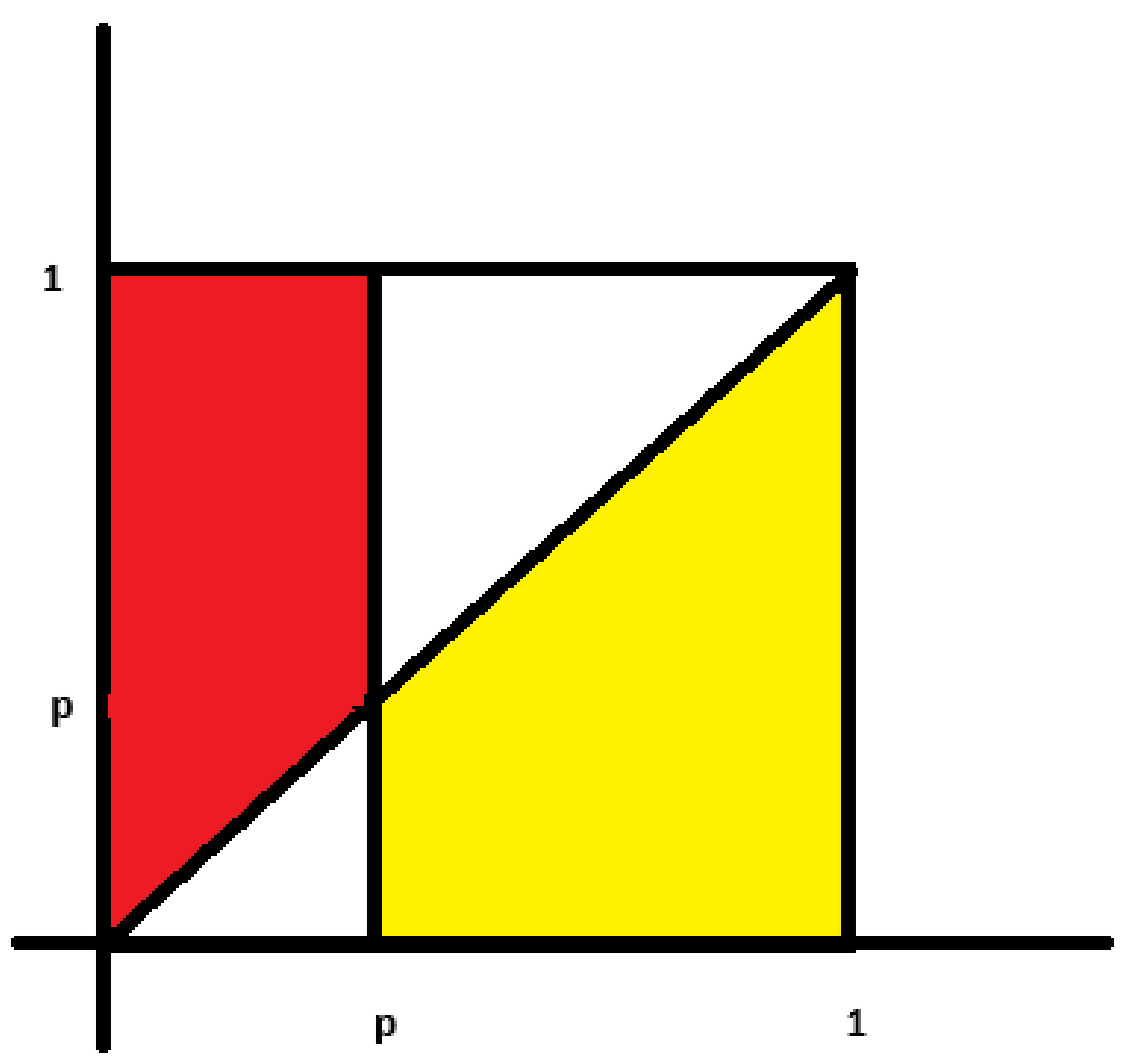

**Fig. 2**

The trapezoid colored red is that portion of the square for which $x < p$ and $y > x$ – so this represents a correct guess. The trapezoid colored yellow is that portion of the square for which $x > p$ and $y < x$ – so this represents a correct guess. So the colored area is $1 - \frac{1}{2} p^2 - \frac{1}{2} (1-p)^2$, and the result follows. █

We now assume the outcomes of Experiment A are uniformly distributed on [a,A] and the outcomes of Experiment B are uniformly distributed on [b.B]. Let the pointer $p \in R$. As previously, we see the result of Experiment A but have no access to Experiment B.

**Lemma 2** – We compute the probabilities of a successful guess in three situations.

(1) [a,A] = [b,B] $P(SG) = F_E(a,A) = \frac{1}{2} + z - z^2$, where $z = \frac{p-a}{A-a}$ if $a < p < A$ and ½ if $p < a$ or $p > A$. The subscript E is chosen to denote the fact that the intervals [a,A] and [b,B] are equal.

(2) $A < b$ $P(SG) = F_L(a,A) = 0$ if $p < a$, $(p-a)/(A-a)$ if $a < p < A$, and 1 if $p > A$. The subscript L is chosen to denote the fact that the interval [a,A] is lower than [b,B].

(3) $B < a$ $P(SG) = F_H(a,A) = 1$ if $x < p < a$. $(A-p)/(A-a)$ if $a < p < A$, and 0 if $p > A$. The subscript H is chosen to denote the fact that the interval [a,A] is higher than [b.B].

**Proof –** (1) This is just a scaling argument from Lemma 1.

(2) If $p < a$, then $p < x$ for all $x \in [a,A]$. Alice therefore guesses that $x > y$, and this guess is always wrong. If $a < p < A$, a correct guess requires Alice to guess $y > x$, and we will do so whenever $x < p$, which occurs with probability $(p-a)/(A-a)$. If $p > A$, then $p > x$ for all $x \in [a,A]$. Alice will therefore guess that $x < y$, and this guess is always right.

(3) If $p < a$, then $p < x$ for all $x \in [a,A]$. We will therefore guess that $x > y$, and this guess is always right. If $a < p < A$, a correct guess requires us to guess $y < x$, and we will do so whenever $x > p$, which occurs with probability $(A-p)/(A-a)$. If $p > A$, then $p > x$ for all $x \in [a,A]$. We will therefore guess that $x < y$, and this guess is always wrong. ∎

It is interesting to realize that if Alice and Bob are conducting the same experiment with outcomes uniformly distributed on some bounded interval, it is not necessary for Alice to know anything about the interval in order to ensure that P(SG) > ½. She need merely choose a pointer which has a non-zero probability of lying in that interval, which could be done simply by selecting a pointer from a countable dense subset with selection probability $1/2^n$ for the $n^{th}$ element of the countable dense subset.

We now use the results of Lemma 2 to obtain P(SG) when the pointer is chosen from a uniform distribution on the interval [p,P] which contains the interval [a,A]. We begin by examining three specific cases. The three cases we are examining are the two cases where the intervals [a,A] and [b,B] do not intersect, and the case in which the interval [a,A] is the same as the interval [b,B].

Each one of these cases has four sub-cases, depending on how the interval [p,P] interfaces with [a,A]. We do not include the cases where $A \leq p$ or when $P \leq a$ because the guesses are always the same – if $A \leq p$, Alice will always guess Experiment B is larger, and when $P \leq a$ she will always guess that Experiment A is larger. We denote these cases by ordering them in order of increasing size. So $F_{LO}(p,a,A,P)$ represents the case where the interval [a,A] does not intersect [b,B] and $A < b$, and $p \leq a \leq A \leq P$.

**Case LO** – $A < b$ – then x is always less than y. We denote the probability of a successful guess in this instance by $F_{LO}$, with further values indicated by the following computations.

There are 4 different allowable orders for a,A,p and P, given that $p < P$ and $a < A$. We denote these by writing down the letters left to right as the values associated with the letters increase.

We first begin by computing $F_{LO}$(p,a,A.P) – the probability of a successful guess when, in addition to A < b, we have p < a < A < P.

x < y and pointer > x - $\int_b^B \int_a^A \frac{P-x}{P-p} \frac{dx}{A-a} \frac{dy}{B-b} = \frac{1}{(P-p)(A-a)} \int_a^A (P-x)dx =$

$$= \frac{1}{(P-p)(A-a)} [P(A-a) - \left(\frac{A^2}{2} - \frac{a^2}{2}\right)]$$

$$\text{So } F_{LO}(p,a,A.P) = \frac{P-(A+a)/2}{P-p}$$

$F_{LO}$(,a,p,P.A) is computed by realizing that when x ϵ [a,p], the pointer is always greater than x, so Alice guesses that Experiment B is larger, and this is always correct. When x ϵ [p,P], half the time the pointer is greater than x and half the time it is less than x, so Alice guesses correctly half the time. When x ϵ [P,p], the pointer is always less than x, so Alice guesses that Experiment A is larger, and this is never correct. As a result,

$$F_{LO}(a,p,P,A) = \frac{p-a}{A-a} + \frac{1}{2}\frac{P-p}{A-a}$$

$F_{LO}$(p,a,P.A) is computed by realizing that when the pointer ϵ [p,a], Alice always guesses that Experiment A is larger and is never right. When the pointer ϵ [a,P], this is a special case of the result in the previous paragraph. As a result,

$$F_{LO}(p,a,P.A) = \frac{1}{2}\frac{(P-a)^2}{(P-p)(A-a)}$$

Finally, $F_{LO}$(a,p,A,P) is computed by realizing that when x ϵ [a,p], the pointer is always greater than x, so Alice guesses that Experiment B is larger, and this is always correct. When x ϵ [p,A], this is a special case of the first formula obtained in the computation of $F_{LO}$. As a result,

$$F_{LO}(a,p,P.A) = \frac{p-a}{A-a} + \frac{A-p}{A-a}\left(\frac{P-\frac{p+A}{2}}{P-p}\right)$$

**Case HI** – B < a – then x is always greater than y

As in Case LO, we first compute $F_{HI}$(p,a,A,P).

x > y and pointer < x - $\int_b^B \int_a^A \frac{x-p}{P-p}\frac{dx}{A-a}\frac{dy}{B-b} = \frac{1}{(P-p)(A-a)}\int_a^A (x-p)dx =$

$$= \frac{1}{(P-p)(A-a)}\left[\left(\frac{A^2}{2} - \frac{a^2}{2}\right) - p(A-a)\right]$$

$$\text{So } F_{HI}(p,a,A,P) = \frac{\frac{A+a}{2} - p}{P-p}$$

$F_{HI}$(,a,p,P.A) is computed by realizing that when x ϵ [a,p], the pointer is always greater than x, so Alice guesses that Experiment B is larger, and this is always wrong. When x ϵ [p,P], half the time the pointer is greater than x and half the time it is less than x, so Alice guesses correctly half the time. When x ϵ [P,p], the pointer is always less than x, so Alice guesses that Experiment A is larger, and this is always correct. As a result,

$$F_{HI}(a,p,P,A) = \frac{1}{2}\frac{P-p}{A-a} + \frac{A-P}{A-a}$$

$F_{HI}$(p,a,P.A) is computed by realizing that when x ϵ [a,,P], this is a special case of the computation of the first formula of $F_{HI}$. When x ϵ [P,A], Alice always guesses Experiment A is larger, and this is always correct . As a result,

$$F_{HI}(p,a,P.A) = \frac{A-P}{A-a} + \frac{P-a}{A-a}\left(\frac{\frac{P+a}{2}-p}{P-p}\right)$$

Finally, $F_{HI}(a,p,A,P)$ is computed by realizing that when the pointer $\epsilon$ [p,A] , Alice is right half the time. When the pointer $\epsilon$ [A,P], Alice guesses Experiment B is larger, and this is always wrong. As a result,

$$F_{HI}(a,p,A,P) = \frac{1}{2}\frac{(A-p)^2}{(P-p)(A-a)}$$

**Case EQ** – a=b, A=B let c = A-a and note that we can assume a = b = 0 by translation-invariance x < y and pointer > x - $\int_0^c \int_0^y \frac{P-x}{P-p}\frac{dx}{c}\frac{dy}{c} = \frac{1}{c^2(P-p)}\int_0^c \int_0^y (P-x)dxdy = \frac{3P-c}{6(P-p)}$ x > y and pointer < x - $\int_0^c \int_y^c \frac{x-p}{P-p}\frac{dx}{c}\frac{dy}{c} = \frac{1}{c^2(P-p)}\int_0^c \int_y^c (x-p)dxdy = \frac{2c-3p}{6(P-p)}$

So

$$F_{EQ}(p,a,A.P) = \frac{1}{2} + \frac{c}{6(P-p)}$$

To compute $F_{EQ}(a,p,P,A)$, we first look at the case where $a=0 \le p \le P \le 1$, Using (2) of Lemma 1 and the formula for the average value of a function, we see that this is

$$\frac{1}{P-p}\int_p^P \left(\frac{1}{2} + x - x^2\right) dx = \frac{1}{2} + \frac{1}{2}(P+p) - \frac{1}{3}(P^2 + Pp + p^2)$$

For $a \le p \le P \le A$, if u = (p-a)/(A-a) and v = ((P-a)/(A-a), the above formula scales to

$$F_{EQ}(a,p,P,A) = \frac{1}{2} + \frac{1}{2}(v+u) - \frac{1}{3}(v^2 + uv + u^2)$$

To compute $F_{EQ}$(a,p,A,P), observe that if the pointer ϵ [p,A], the probability of a correct guess has the form $F_{EQ}$(a,p,P,A) where A = P – and we have just computed this. If the pointer ϵ [A,P], then Alice will always guess that Experiment B is larger, and she will be correct half the time.

We let u = (p-a)/(A-a). Then

$$F_{EQ}(a,p,A,P) = \frac{1}{2}\frac{P-A}{P-p} + \frac{A-p}{P-p}\left(\frac{2}{3} + \frac{1}{6}u - \frac{1}{3}u^2\right)$$

Finally, to compute $F_{EQ}$(p,a,P,A), observe that if the pointer ϵ [p,a], then Alice will always guess that Experiment A is larger, and she will be correct half the time. If the pointer ϵ [a,P], the probability of a correct guess has the form $F_{EQ}$(a,p,P,A) where a = p – and we have just computed this.

We let v = (P-a)/(A-a). Then

**When Experiments A and B Have Overlapping Uniform Probability**

We can use the material in the preceding section to obtain expressions for the probability of a successful guess when the distributions governing the outcomes of Experiments A and B differ, but both are uniform. The following theorem tabulates the probabilities P(SG) in the case that the intervals that represent the ranges for the outcomes overlap

The formulas below are derived from the functions defined above by simply decomposing overlapping intervals into three non-intersecting intervals.

**Case 1 – a < b < A < B**

In this case, Alice's result either belongs to the interval [a,b] or the interval [b,A]. If Alice's result belongs to the interval [a,b], then [a,b] acts as the predictor experiment for $F_{LO}$ on [b,B].. If Alice's result belongs to the interval [b,A], we break up the interval [b,B] into the union of the intervals [b,A] and [A,B]; then [b,A] acts as the predictor experiment for $F_{EQ}$ on [b,A] and as the predictor experiment for $F_{LO}$ on [A,B].

The formulas for the other cases are arrived at in a similar fashion. Each of these overlapping cases allows for strict prediction as long as Experiment A has known values for a and A, as we can then choose P=A and p=a. The proof of the easiest case (Case 2) is given below, the rest are left to the reader.

Probability of successful guess P(SG) =

$$\frac{b-a}{A-a}F_{LO}(a,b)+\frac{A-b}{A-a}\left(\frac{A-b}{B-b}F_{EQ}(b,A)+\frac{B-A}{B-b}F_{LO}(b,A)\right)$$

The notation used here suppresses the role of [p,P] – but it is the relationship between [p,P] and [a,A] which dictates which subcase of the formulas for $F_{LO}$ and $F_{EQ}$ are used. For example, suppose that the subcase being used has $p < a < P < A$. Then the first variable used for $F_{LO}$ and $F_{EQ}$ in the above equation describe which variable should substitute for a in $F_{LO}(p,a,P,A)$ and $F_{EQ}(p,a,P,A)$, and the second which variable should substitute for A. Thus, if [p,P] and [a,A] satisfy $p < a < P < A$, the above expression indicates that (left to right)

(1) a substitutes for a and b substitutes for A in $F_{LO}(p,a,P,A)$

(2) b substitutes for a and A substitutes for A in $F_{EQ}(p,a,P,A)$

(3) (2) b substitutes for a and A substitutes for A in $F_{LO}(p,a,P,A)$

This convention is followed in the remaining subcases.

**Case 2 – b < a < A < B**

Probability of successful guess P(SG) =

$$\frac{a-b}{B-b}F_{HI}(a,A)+\frac{A-a}{B-b}F_{EQ}(a,A)+\frac{B-A}{B-b}F_{LO}(a,A))$$

The proof that this results in a strict prediction as long as P = A and p = a follows.

P(SG) =

$$\frac{a-b}{B-b}\left(\frac{\frac{A+a}{2}-p}{P-p}\right)+\frac{A-a}{B-b}\left(\frac{1}{2}+\frac{A-a}{6(P-p)}\right)+\frac{B-A}{B-b}\left(\frac{P-\frac{A+a}{2}}{P-p}\right)=$$

$$\frac{a-b}{B-b}\left(\frac{\frac{A+a}{2}-a}{A-a}\right)+\frac{A-a}{B-b}\left(\frac{1}{2}+\frac{A-a}{6(A-a)}\right)+\frac{B-A}{B-b}\left(\frac{A-\frac{A+a}{2}}{A-a}\right)=$$

$$\frac{a-b}{B-b}\left(\frac{1}{2}\right)+\frac{A-a}{B-b}\left(\frac{2}{3}\right)+\frac{B-A}{B-b}\left(\frac{1}{2}\right)=$$

$$\frac{a-b}{B-b}\left(\frac{1}{2}\right)+\frac{A-a}{B-b}\left(\frac{1}{2}\right)+\frac{B-A}{B-b}\left(\frac{1}{2}\right)>\frac{1}{2}$$

**Case 3 – $a < b < B < A$**

Probability of successful guess P(SG) =

$$\frac{b-a}{A-a}F_{LO}(a,b)+\frac{B-b}{A-a}F_{EQ}(b,B)+\frac{A-B}{A-a}F_{HI}(B,A))$$

**Case 4 – $b < a < B < A$**

Probability of successful guess P(SG) =

$$\frac{A-B}{A-a}F_{HI}(B,A)+\frac{B-a}{A-a}(\frac{B-a}{B-b}F_{EQ}(a,B)+\frac{a-b}{B-b}F_{HI}(a,B))$$

We note that analogous formulas – which are somewhat messier – can be derived from Lemma 2 when the pointer is a single point.

**Ex. 7 –** Suppose that Alice is conducting Experiment A, the uniform distribution on [2,3], and uses the uniform distribution on [2,3] as the pointer. If Bob is conducting Experiment B, the uniform distribution on [1,4], the Case 2 formula shows that P(SG) = 5/9. Now suppose the situations are reversed and Alice is conducting Experiment A, the uniform distribution on [1,4] with the uniform distribution on [1,4] as the pointer. If Bob is conducting Experiment B, the uniform distribution on [2,3], the Case 3 formula shows that P(SG) ≈ 0.74.

**Section IV – Simulations**

The literature leaves one with the impression that having an equal choice of which of the two pieces of paper to examine is central to Cover's guessing game, and the results of this paper show that, whenever Bob's numbers are obtained via an experiment based on a probability distribution, it is not. Initially the author found these results disquieting, and decided to conduct simulations. The programs were written in BASIC, but the essence of the simulation was that Experiment A was performed and the guessing strategy implemented before any action was taken on Experiment B. Computer instructions were executed sequentially as indicated above. All results in this paper were confirmed by chi-square goodness-of-fit tests.

The results of the simulations extend beyond the results presented here. What those results showed was that in many instances, all that needs to be known about Experiment A is the range of the outcomes; if one invokes a uniform (or normal with midpoint of the range as the mean and a standard deviation approximately a quarter of the range) pointer on that range, Experiment A predicts the result of a wide number of Experiment Bs (in the sense of this paper) with a probability greater than ½; in many cases significantly greater than 1/2 .

These simulations suggest that the results obtained in this paper can be substantially generalized and extended.